# Ethical and sustainable mathematics is localised: why global paradigms fail and culturally-situated practices are essential

Dennis Müller[1] & Maurice Chiodo[2]

7 October 2025

## Abstract

This paper identifies several different interconnected challenges preventing the move towards more ethical and sustainable mathematics education: the entrenched belief in mathematical neutrality, the difficulty of simultaneously reforming mathematics and its pedagogy, the gap between academic theory and classroom practice, and the need for epistemic decolonisation. In this context, we look at both bottom-up and top-down approaches, and argue that globalised frameworks such as the United Nations' Sustainable Development Goals are insufficient for this transformation, and that ethical and sustainable forms of mathematics ought not to be built using these as their (philosophical) foundation. These frameworks are often rooted in a Western-centric development paradigm that can perpetuate colonial hierarchies and fails to resolve inherent conflicts between economic growth and ecological integrity. As an alternative, this paper advocates for embracing localised, culturally-situated mathematical practices. Using the Ethics in Mathematics Project as a case study within a Western, Global North institution, this paper illustrates a critical-pragmatic, multi-level strategy for fostering ethical consciousness within a specific research community, and shows how this may be achieved in otherwise adversarial circumstances.

[1] Institute of Mathematics Education, University of Cologne, Cologne, Germany. dennis.mueller@uni-koeln.de **[corresponding author]**
[2] Centre for the Study of Existential Risk, University of Cambridge, Cambridge, United Kingdom. mcc56@cam.ac.uk; ORCID: https://orcid.org/0009-0006-3479-7822

# Introduction

As the societal impact of mathematics becomes increasingly visible (e.g., Chiodo & Clifton, 2019; O'Neil, 2017; Skovsmose, 2024), researchers are increasingly discussing the mathematicians' and educators' ethical and sustainable responsibilities (Budikusuma et al., 2024; Müller, 2025a). In particular, the relationship between mathematics education and climate change presents a layer of complexity that demands urgent attention (Barwell, 2013; Skovsmose, 2021). Sustainability concerns are "wickedly" complex, which means that time is running out to find comprehensive solutions, those who initially cause the problem are also tasked with creating such solutions, there is no central authority strong enough to enforce a particularly sustainable response, and policies tend to discount the future for irrational reasons (Levin et al., 2012). Such "wicked" problems resist simplicity and often spark disagreement in society (Hulme, 2009), demanding reconsiderations of the mathematical curriculum and pedagogy to work towards sustainable futures.

When these topics are introduced into classrooms, their "wicked" features accompany them, pushing teachers and students to abandon perceptions of mathematical neutrality (e.g., Fred et al., 2024a, 2024b; Steffensen et al., 2018, 2023). This abandonment occurs because wicked problems inherently involve conflicting values and stakeholders, meaning that any mathematical problem, model, and solution will inevitably prioritise certain interests over others, thus exposing the inherent political nature of mathematical choices and artefacts (Marichal, 2025; Müller & Chiodo, 2023). Consequently, traditional approaches to mathematics that emphasise the subject's certainty, precision, and objectivity appear increasingly inadequate for preparing students to engage with complex socio-ecological realities. Valero (2023) argues that mathematics, as a central project of modernity, emphasises control and exploitation of nature, casting it merely as an object to be studied and handled. Similarly, Barwell (2024) suggests that the mathematical consciousness cultivated in traditional schooling may contribute directly to today's crises. An example of this is given by the oil pipeline problem described by Shulman (2002, p. 124), which frames resource extraction as a matter of profit-making:

> An offshore oil well is located at point A, which is 13 kilometers from the nearest point Q on a straight shoreline. The oil is to be piped from A to a terminal at a point T on the shoreline by piping it straight under water to a point P on the shoreline between Q and T and then to T by a pipe along the shoreline. Suppose that the distance QT is 10 kilometers, that it costs $90,000 per kilometer to lay underwater pipe, and that it costs $60,000 per kilometer to lay the pipe along the shoreline. What should the distance be from P to T in order to minimize the cost of laying the pipe?

This framing implicitly validates the extraction process and prioritises economic efficiency while rendering ecological and social consequences invisible. Its author exemplifies how mathematical training can normalise unsustainable practices by obscuring the normative assumptions embedded within seemingly neutral problems. The exercise's task aligns with broader critiques that educational institutions often perpetuate the cultural patterns that drive the ecological crisis (Bowers, 2017). Shulman (2002, p. 124) and Chiodo et al. (2025, pp. 61 - 62) have thus suggested to augment the question by asking what other factors need to be considered, in order to "get students to realise that there are at least SOME factors, and that optimising over (naive) cost, with no other considerations, is not at all a good way to answer the question" (ibid., p. 62).

On a broader level, the very imposition of universal curricula bringing such language and mathematical consciousness into the world's classrooms can be understood as authoritarian acts of control that can be harmful for students (Ernest, 2024b) and society (Chiodo & Müller, 2025). In response to such shortcomings, different visions for mathematics education are emerging, in which the curriculum ought to prioritise critical judgement. Among other methods, they may focus on "normative modelling" (Pohlkamp, 2022) or "embedded ethics" (Müller & Chiodo, 2023), giving students the chance to (continuously) face ethical and sustainability concerns in the mathematics exercises they encounter. Such mathematical practice and teaching involves not just solving a prescribed problem, but also understanding its underlying construction. For instance, analysing carbon footprint calculations can reveal the implicit assumptions and competing interests that shape their framing (Weiss & Kaenders, 2023).

As described by Coles et al. (2024), the socio-ecological entanglement of today's crises calls for transcending the artificial separation between the social and ecological dimensions and to recognise their interconnection and mutual constitution. In other words, ethics and sustainability should be considered and taught within mathematics, rather than separately (Chiodo & Bursill-Hall, 2019). Mathematics that engages critically with sustainability and ethics can empower students to analyse complex socio-ecological problems (Meyer et al., 2025), interrogate power structures (Ernest, 2002, 2019), and help to imagine alternative futures (Makramalla et al., 2025). To achieve this within mathematics, one sees that mathematicians, teachers, and students need to hone their ability to navigate value-laden decisions rather than pursuing ostensibly neutral calculations and proofs (Rycroft-Smith et al., 2024). Such perspectives on ethical and sustainable mathematics align with broader critiques arguing that mathematics, as a human practice, is neither objective, nor neutral, nor value-free (e.g., Chiodo & Müller, 2024b; Ernest, 2007, 2020, 2016, 2024a; Mellin-Olsen, 2002).

Drawing on recent developments in critical mathematics education, ethnomathematics, and sustainability studies, we will begin by discussing some of the interconnected systemic challenges the subject is currently facing regarding ethics and sustainability. We will argue that sustainable mathematics demands a fundamental reconsideration of the philosophical foundations of mathematics, the simultaneous transformation of both mathematical practice and pedagogy, the bridging of persistent gaps between theoretical frameworks and classroom realities, and a commitment to epistemic decolonisation that challenges Western-centric conceptions of mathematical knowledge.

We then critically examine whether global frameworks, such as the United Nations' (2015) Sustainable Development Goals, provide adequate guidance for this transformation, ultimately arguing for the necessity of localised, culturally-situated approaches to mathematical practice. Through the example of the Ethics in Mathematics Project (2025), we illustrate both the potential and limitations of local, critical-pragmatic interventions within specific institutional contexts, emphasising that meaningful change requires not merely modifying existing structures but fundamentally reimagining the relationships between mathematics, culture, power, and knowledge production, and the teaching thereof.

In the following, we will say something is "localised" when it is deliberately adapted for a specific local community. The Ethics in Mathematics Project (EiMP) is an example of a localised approach because its frameworks were strategically tailored to the specific culture of research mathematicians, enabling new forms of socio-cultural locality to emerge in which shared interests on ethics and sustainability can develop (for a further discussion of localisation, we refer to Cox & Mair, 1991). Our usage of "epistemic decolonisation" follows (Leiviskä, 2023, p. 227) who defined it as "challeng[ing] and transform[ing] the epistemological foundations involved in the production and transmission of knowledge". Our work is, thus, concerned with tacit forms of mathematical practice and their teaching. This represents what Malik and Neuweg (2025, pp. 43-44) call a "weak reading" of decolonisation, one not primarily focused on overcoming oppression.

# On some interconnected challenges

Implementing any vision of ethical and sustainable mathematics is not straightforward. While there is evidently an emerging consensus among education researchers who see the urgency to respond to ecological collapse (Oikonomou, 2025), fundamental questions remain about the practical nature and adequacy of the proposed responses. A central debate revolves around the depth of this integration: are we only going to teach *about* ethics and sustainability, advocate *for* specific forms of ethics and sustainability, or rebuild mathematics *as* ethical and sustainable (Renert, 2011)? In other words, the critical questions for research into ethical and sustainable mathematics education have moved away from only asking *whether* one should integrate ecological and justice concerns into mathematics education, to the more difficult question of *how* it can be done, and how (or if) consensus can be established (cf. Barwell & Hauge, 2021; Chiodo & Müller, 2024a).

## Practitioner resistance and the myth of neutrality

Both teachers (cf. Abtahi et al., 2017; Gray & Bryce, 2006, Ernest 2016) and mathematicians (cf. Chiodo & Bursill-Hall, 2018, 2019) can demonstrate restraint about the political nature of ethics and sustainability in mathematics education. As Steffensen et al. (2018, p. 233) argue: "Political guidelines and curricula can put teachers in a challenging and risky situation. Mathematics has traditionally been regarded as a neutral subject with little controversy. This can make it extra challenging to include wicked problems like climate changes [sic!] in mathematics lessons." And students may develop resistance to such teaching if their educators themselves are resistant (Chiodo & Bursill-Hall, 2019). Such behaviour can have various origins, but it likely includes balancing philosophical disagreements about the nature of mathematics, concerns about politicising the classroom, the need to learn challenging new skills, or simply the overwhelming demands already placed on educators (cf. Chiodo & Müller, 2024a; Lestari et al., 2024; Su et al., 2022, 2023).

Indeed, reviews suggest that the training of teachers and mathematicians often only contains superficial lessons on ethics and sustainability concerns, indicating a need to build deeper competencies among the subject's practitioners (cf. Su et al., 2023; Müller et al., 2022). However, among mathematics teachers, a lack of knowledge about sustainability is also reported in circumstances where the educators value the connection between mathematics and sustainability (Alsina & Vásquez, 2025). This lack of competency may be exacerbated by structural constraints, such as accreditation requirements prioritising canonical content over critical pedagogy, and the professional risks educators may face when addressing controversial topics without explicit administrative support. To meet these challenges, new teaching resources on ethics and sustainability, including university-level textbooks (e.g., Roe et al., 2018; Chiodo et al., forthcoming) and school resources (e.g., Schwarz et al., 2024; Mathematics Hub, n.d.), are being developed.

In particular, the resistance from deep-seated purist and absolutist beliefs about mathematics presents a significant barrier to dismantling myths of neutrality (e.g., Burton 1994; Ernest, 2021a; Wagner, 2023). Yet, while there is hope that such myths are being challenged from within the discipline (Tractenberg et al., 2024; Buell et al., 2022), the disconnect between research and practice still leads to uneven results when the radical critiques emerging from academic discourse encounter practical constraints that educators face in their classrooms, resulting in what Li and Tsai (2022) describe as "patchwork implementations" — fragmented approaches which fail to realise the transformative potential of critical mathematics education.

## Shaping the future: existential sustainability and pedagogical visions

Beyond these practical hurdles, the ongoing debates touch on the fundamental role that mathematics plays in shaping our perception of the future. In today's capitalist societies, mathematics is fundamental to how we imagine the future, and enables us to make decisions in uncertain times (Beckert, 2016; Beckert & Bronk, 2018). Modelling and optimisation, in particular, shape our understanding of how the climate and other sustainability crises might evolve (e.g., Barwell, 2013; Geiger, 2024; Skovsmose, 2021; Marichal, 2025). Teaching nothing, or only very little, about ethics and sustainability thus limits the students' abilities to imagine and think in alternatives. The metacognitive ability to understand both how to use mathematics, and how mathematics shapes thought, language, and action, becomes essential for meaningfully engaging with ethics and sustainability. Thus, following Lakatos's (1978, p. 20) observation that we build our own (epistemic) prisons, but also possess the critical capacities to destroy them, D'Ambrosio (2015, p. 24) proposes that mathematics education must encompass the critique of, and connection with, our everyday existential experiences.

Such attempts introduce what Müller (2025b) refers to as "existential sustainability" into the classroom. That is, a paradigm of mathematical practice, learning, and teaching that values mathematical questions, sustainability concerns, and the existential well-being of students (and other humans and species) from a nuanced and holistic perspective. Besides imagination, hope has become a central feature of critical education (Gottesman, 2016).

Such approaches move beyond simply applying mathematics to sustainability problems; they address how the very act of practicing mathematics shapes students' sense of self, their agency, and their personal relationship with the world. Increasingly, one can observe that the subject calls for a curriculum design that genuinely sees students as multi-dimensional beings and aims to prepare them for alternative, sustainable futures, rather than being stuck in specific concrete normative, or epistemic visions (Makramalla et al., 2025).

Winter (1996, p. 35) explains that a general mathematics education ought to provide three foundational experiences to students: perceive and understand socio-planetary and cultural concerns, experience mathematics for its own sake, and build problem-solving skills allowing students to apply mathematical heuristics in life. Regarding ethical and sustainable general mathematics education, the heart of this challenge then becomes the need to foster a dialogic relationship with the living world (Barwell et al., 2022), through which students live with nature rather than seeing it as an object of mathematical desires (cf. Gutiérrez, 2017; Müller, 2025b; Chiodo & Müller, 2025; Müller & Chiodo 2023). All of this reveals a problem that operates on multiple levels: micro questions about the correct exercise design and teaching choices are entangled with meso-level challenges of institutional constraints, and macro issues of the globally prevalent philosophical beliefs surrounding mathematical neutrality. This raises crucial questions about how pragmatic and critical perspectives can come together in mutually beneficial ways in the classroom (Müller, 2025b). The challenge is not merely one of knowledge transfer, but also of creating the conditions that enable teachers to translate sophisticated theoretical insights into meaningful learning experiences within existing institutional structures and constraints, to form an active citizenship education (e.g., Maass et al., 2019; Buell & Shulman, 2019).

## Bridging the theory-practice gap: bottom-up approaches

Some suggest that bottom-up approaches are suitable here, focusing on ethically valued content that connects with students' lived experiences to build agency (e.g., Buell & Piercey, 2024). This insight is supported by early research indicating that students' experience of self-efficacy positively influences their willingness to engage with complex sustainability problems, suggesting that tasks should be designed to promote confidence and focus on achievable successes (Meyer et al., 2025). While mathematical exercises are a commonly discussed topic, bottom-up approaches also focus on classroom culture. For example, Wilhelm (2024) argues for implementing sustainable classroom cultures through mindful teaching and Coles (2023) calls for "dialogic ethics" as part of teaching. But, while feasible within local settings, Müller (2025a) notes that many teachers are not taught how to do this, thus requiring bravery and resilience to learn such skills and to deploy them.

Chiodo & Bursill-Hall (2019) note that such bottom-up teaching requires the alignment of three institutional concerns: the integration of ethics directly into mathematical exercises, teaching students about ethics outside of their normal mathematics classes, and gathering the support of other teaching staff. Buell & Shulman (2019) argue that socially-just pedagogy has three components: socially-just curriculum, methodology, and education practices. Regarding ethical and sustainable teaching, the critical question which lies behind both Chiodo & Bursill-Hall's (2019) and Buell and Shulman's (2019) perspective on ethical teaching, is how such cultural shifts can be fostered within existing educational systems that often resist fundamental change, particularly when these changes challenge established power dynamics and assessment regimes.

To address this challenge, Rycroft-Smith et al. (2024, p. 377) proposed focusing on "usable frameworks", suggesting a level-based guideline for educating teachers, slowly moving them from "obstructing efforts to address ethics in their mathematics classroom" to "taking a seat at the tables of power" where they can advocate for (local) change and use their "ethical awareness to object to bad data in education systems; questioning and contributing to policies; exposing mathwashing; auditing and challenging with mathematics; confronting gatekeeping".

## The imperative of epistemic decolonisation

Zooming out from the classroom and considering global dynamics, there is always the question of whose ethics we are following (Dubbs, 2020). The focus on "usability" by Rycroft-Smith et al. (2024) must be critically studied with regards to the question of whether such a utilitarian focus accidentally risks perpetuating existing harmful path dependencies. As we discuss in the next two sections of the paper, localised, bottom-up approaches can be extremely well suited to break with colonial legacies of power. Engaging in genuine epistemic decolonisation efforts represents a complex terrain for educators and researchers to navigate. While ethnomathematics manages to offer valuable suggestions on incorporating cultural and linguistic resources of non-Western communities into mathematics education (Rosa & Orey, 2016), scholars have also identified significant risks related to forms of superficial appropriation that may ultimately reinforce rather than challenge the dominance of conventional Western mathematics (e.g., Abtahi, 2022; Chronaki & Lazaridou, 2023; Swanson & Roux, 2025). The central tension lies in how these cultural resources are integrated: when used merely as motivational contexts for teaching conventional mathematics, the underlying power dynamic remains unchanged, and the epistemic richness of these traditions is marginalised rather than centered.

Thus, both mathematics and mathematics education must grapple with how to authentically challenge Eurocentric definitions of mathematical practice, teaching and learning that propose mathematical universality and neutrality while avoiding pitfalls of tokenism or cultural extraction. This requires not merely adding diverse content to existing teaching frameworks or focusing on "usable" ideas, but more fundamentally to reconsider what counts as mathematical knowledge and who has the authority to make such determinations (e.g., Müller, 2024). For example, recent scholarship has begun to question whether current equity frameworks in mathematics education may be fundamentally flawed due to their underlying assumptions. Bullock (2023) argues that the very logic structuring both school mathematics and equity research may inadvertently undermine their stated aims by implicitly accepting an "axiom of racialised deviance", perpetuating colonial epistemologies. Ultimately, such arguments force mathematics and its educational studies to confront an uncomfortable possibility: that incremental reforms within existing structures can be insufficient or even counterproductive if they do not also address larger philosophical and foundational concerns related to mathematics itself. So what can be said about top-down approaches, such as those that build on the Sustainable Development Goals?

## Top-down approaches: are the Sustainable Development Goals good enough to achieve this?

In a comprehensive meta-analysis of over 3,000 scientific studies, Biermann et al. (2022) assess the political impact of the United Nations' Sustainable Development Goals (SDGs) and conclude that their influence has been limited and largely "discursive". While the SDGs have successfully provided a common language for sustainable development in international relations and national politics, their introduction has not led to widespread, tangible change. Biermann et al. find little evidence that the goals have had a direct, transformative impact on politics in terms of substantive legislative action, resource reallocation or the realignment of existing institutions with sustainable futures. Instead, governments occasionally implement the SDGs selectively, choosing goals that align with pre-existing agendas, while sub-national actors, such as cities and private institutions, tend to exhibit more progressive adoption. This fragmented implementation is compounded by a struggle to achieve policy coherence and consistent integration of the SDGs. Despite some advances, progress is frequently stalled by bureaucratic inertia, short-term political interests, and a diminishing sense of ownership over the goals.

Furthermore, Biermann et al. find that the SDGs have fallen short in their substantive commitments to social justice. The core principle of "leaving no one behind" has not translated into the envisioned reduction in social inequality, and, in some instances, the SDGs have even been used to legitimise existing marginalisation. They further point to some inherent contradictions within the SDGs, such as the conflict between universal economic growth and environmental protection targets, giving nations the excuse to prioritise socioeconomic objectives over ecological concerns.

Thus, it appears that an insurmountable conflict between the SDGs and holistic forms of ethical and sustainable mathematics education lies at the level of epistemology. The SDGs form the latest iteration of development paradigms that threaten to perpetuate specific, Western-centric ways of knowing and socio-economic organisation;[3] a perspective which is irreconcilable with the call for epistemic decolonisation, the diversity of thought, and the fundamental critique of mathematics presented in the previous section. Post-colonial scholarship has criticised the very concept of "development" as a dividing act, producing global hierarchies and reproducing older colonial relations: those of the Global South are in some sense less "developed", and reliant on knowledge and aid from the Global North (Vogt, 2022).

Krauss et al. (2022) are particularly sceptical of SDGs 8 (decent work and economic growth), 9 (industry, innovation, and infrastructure), 12 (responsible consumption and production), 13 (climate action), and 15 (life on land), arguing that the United Nations may have continued in the footsteps of what Quijano (1992, 2000) called "coloniality of power". These footsteps are not just economic but also linguistic, as Western assumptions about progress and individualism are embedded in the very language of development, undermining other cultural ways of knowing (e.g., Bowers, 2017). Crucially, the SDGs implicitly rely on the perceived universality and neutrality of mathematics to measure, monitor, and define progress, thereby reinforcing precisely those assumptions that truly sustainable mathematics seeks to dismantle.

---

[3] For example, Zai (2021) writes that "the SDGs follow entirely in the tradition of development policy, which was established in the context of the Cold War and the decolonization of the USA and its allies with the aim of preventing those countries that were becoming independent from defecting to the socialist camp. This was done by promising them that they could become, with the support of the West, prosperous or "developed" countries within a capitalist world economic order."

Regarding discourses of education scholars, Anuradha (2024) finds that epistemic injustices within sustainable development have become a central concern, leading them to argue that social justice and epistemic injustice cannot be solved independently, but rather must always be considered simultaneously. Practical pedagogical approaches that address this include teaching for spatial justice, which uses mapping tools and geographic information systems to empower students to analyse and critique inequities within their own communities (e.g., Rubel et al., 2017). And indeed, where we see the SDGs perform better is at the level of task and activity design (e.g., Meyer et al., 2025; Schwarz et al., 2024). Some teacher training programs have successfully used the SDGs at the pedagogical level as a framework for designing interdisciplinary STEAM activities that enhance educators' sustainability competencies (e.g., Alsina & Silva-Hormazábal, 2023), suggesting their utility as a practical tool, even if they cannot lead to a deeper answer to the question of "What is sustainable mathematics, and how can it be philosophically and critically conceptualised?"

Therefore, while the SDGs might serve as a common language for introducing sustainability concepts, they should not be mistaken for the endpoint of curriculum transformation, nor as new philosophical foundations for ethical and sustainable mathematics. The need to reconstruct the philosophy of mathematics and its education on ethical and sustainable foundations paradoxically represents the SDGs' antithesis, as it constitutes a project of epistemic decolonisation and epistemic justice. By interrogating what counts as ethical and sustainable mathematical knowledge and practice, and who has the authority to make such determinations, it necessitates a departure from Western paradigms of mathematical competence towards a set of alternative epistemologies, recognising and valorising diverse mathematical traditions and knowledge systems as autonomous and legitimate, possessing inherent validity and power. In short, it aims to overcome mathematics as what Bishop (1990, p. 51) termed the Western world's "secret weapon of cultural imperialism".

Thus, while the SDGs can guide individual teachers, curricula, and the design of activities, the deeper philosophical incompatibility appears rather absolute: initiating a project of ethical and sustainable mathematics on new foundations can only occur as a form of epistemic liberation outside the colonialities of power of modern development programmes. Otherwise, the approach effectively re-colonises mathematics even before genuine efforts have taken place. Unlike admonishing the Lakatosian prison, it would (re-)create an existing cognitive cage. But if a globalised approach building on the SDGs is unlikely to succeed, what can be done instead?

# On the need for localised mathematical practices

Ethnomathematical research has shown that mathematics does not constitute a singular, monolithic discipline but rather a constellation of culturally-situated practices and disciplinary subcultures (Borba, 1990; Bowers, 2008; D'Ambrosio, 2006; Rosa et al., 2016; Rosa & Orey, 2016). Different sociocultural contexts generate mathematical practices calibrated to their particular constraints and axiological commitments, turning mathematics into a product of their local cultures. Thus, Rowlands & Carson (2002) argue that only (academic) mathematics informed by ethnomathematical insights can develop its full potential. Such a recognition that mathematical knowledge production is inherently anthropological and culturally embedded creates the space to reconceptualise its practice through frameworks emphasising localisation and local (Indigenous) philosophies.

The power of such localised approaches is not merely theoretical; it is demonstrated in a growing body of work where mathematics education is re-rooted in the socio-ecological, ethical, and cultural concerns of specific communities. These directly counter some of the weak spots of global development frameworks. Where the SDGs impose the external Western epistemological metrics associated with "development", these initiatives begin by grounding mathematics in the immediate concerns and epistemologies of different communities. They do not seek to integrate a community's local knowledge into the dominant mathematical paradigm merely as context; rather, they want to bring this knowledge forward to establish new critical foundations for mathematical inquiry, i.e., to empower communities to define problems and develop solutions on their own terms, and to counter the "coloniality of power".

For instance, projects in rural Mexico have brought together teachers, scientists, and activists to design interdisciplinary activities using socio-critical mathematical modelling to address immediate local issues like industrial river pollution (Solares-Rojas et al., 2022). Similarly, a pedagogic workshop in a rural Greek village has centered on local knowledge and "the commons" to "re/make" space for learning and subvert the erasure of local epistemic practices (Chronaki & Lazaridou, 2023). In an Indigenous (Canadian) teachers' training course, focussing on local issues helped Abtahi (2022) to overcome communication barriers and establish a sustainable classroom culture. Other initiatives have focused on valuing students' existing mathematical "funds of knowledge" by documenting the mathematics present in the homes and communities of Pacific learners (Hunter, 2022), have explored place- and land-based pedagogies as a direct means of decolonising mathematics education (Nicol et al., 2020), or created specific projects to teach socially-just mathematics to Black communities in the southern parts of the United States (Moses & Cobb, 2002).

These examples illustrate that localised practices can take many forms but share a common goal: to turn mathematics into a meaningful tool for understanding and shaping the world of local communities, rather than decontextualised universal procedures of a "Nearly Universal and Conventional (NUC) mathematics" (Barton, 2008). Instead of presenting students with out-of-context mathematics, they connect students with their "out-of-school worldview", effectively conserving and valuing a community's cultural heritage (François et al., 2018). Crucially, the ethnomathematical understanding of "local" extends beyond geography and ethnicity; it applies to any community defined by its shared practices, values, and constraints, including distinct (academic) disciplinary subcultures. Therefore, a "local" context can be a rural village in Mexico, an Indigenous community in Canada, or even the specific academic culture of a research university in the Global North. This expanded definition is essential, as it frames the case study of the next section. It suggests that even within the institutions that produce and export dominant mathematical paradigms, meaningful change requires interventions that are carefully tailored to the local epistemological and normative commitments as well as to the institutional realities of that specific community.

In other words, "all mathematics is ethnomathematics as a development and result of societies and individuals making part of and expressing through cultural groups of any ethnic composition" (Pareyon, 2022, p. 9). Treating today's research mathematics and its teaching not as a universal, objective, a-cultural truth opens up new venues to incorporate ethics and sustainability. By recognising the mathematics taught in research universities and the community of research mathematicians as a distinct socio-cultural group with its own language, jargon, and symbols (Rosa & Orey, 2011) and codes of behaviour (Müller et al., 2022), one can also analyse and change the ethics behind the NUC mathematics in such places. Therefore, even mathematics in institutions of the Global North can be analysed through the lens of localised practices, once one moves past the naive misconception that such ideas are only for Indigenous communities and the Global South. And in return, such a local focus within seemingly universal research mathematics allows for the development of specific tools and methods aimed at breaking down the existing assumptions and beliefs about the neutrality of mathematics.

However, advocating for localised practices is not without risks. A potential problem lies in its danger to overly romanticise local contexts, which may overlook internal power dynamics or reinforce parochialism if taken on uncritically, and in the fact that local efforts may lack sufficient resources and knowledge required for ethics and sustainability. Furthermore, localised approaches can be difficult to scale and can struggle to address transnational problems. Mathematics as a global field of science and applications, and as a local field of education, will always exhibit such tensions. Thus, critical scrutiny is warranted when mathematicians claim that the criteria for judgement have to come from within the mathematical "profession" (e.g., Borovik, 2023). Such perspectives rely on implicit assumptions that the mathematical community is able to establish neutral, acultural, and universally applicable perspectives -– the very aspect which decades of ethnomathematical, historical, and philosophical scholarship have tried to dissect to showcase a culturally-situated, value-laden discipline. Thus, in what follows, and when analysing the Ethics in Mathematics Project (EiMP), we have to be careful to make its local nature visible, and to carefully trace some of its important external influences.

## The Ethics in Mathematics Project

In the following, we analyse how the EiMP, as one such project focused on university mathematics, has tried to bring ethics and sustainability to mathematicians. Acknowledging the inadequacy of universalist neutrality paradigms within mainstream mathematics, the EiMP developed a sequence of frameworks designed to cultivate ethical consciousness in professional mathematicians and mathematics students. Crucially, one notes that the EiMP was preceded by the Cambridge University Ethics in Mathematics *Society* (CUEiMS, n.d.); a formally-recognised student society within the University of Cambridge, set up in 2016 (two years prior to the formation of the EiMP) to initiate and support the development and promotion of the very notion that there exist ethical issues in mathematics. It was this very clear bottom-up approach, of a student society formed *by* students, *for* students, which gave a sufficient platform and support for the more academically-based EiMP to form. The EiMP itself was founded by the authors with the help of a local historian of mathematics.

Such grassroots origins enabled the EiMP to form and evolve with its audience of mathematics students in mind. As such, it did not produce an understanding of ethics in mathematics "for the sake of it", but rather with a view to resonating with, and teaching, the very students who would one day go on to use (and teach) such knowledge and understanding. Over the years, as the EiMP grew and expanded, it retained the material and moral support of the students from CUEiMS, as evidenced by their continued hosting since 2016 of the (non-exmainable, optional) annual seminar series on "Ethics for the Working Mathematician" given by members of EiMP (see CUEiMS, 2025). Critically, however, the society's guest speaker series, the organisation of two interdisciplinary workshops, and the early collaboration with history of science brought in external perspectives about the localised nature of the project's approach to ethics (Müller, 2024). In particular, as Müller (2024) explains, it showed how even within the academic landscape of the Global North, intellectual disagreements can happen when researchers from different institutions and research backgrounds meet to discuss ethics and social justice in mathematics.

While the EiMP arguably does not represent the ultimate vision of a fully decolonised ethical and sustainable mathematics (cf. Ernest, 2021b), it is presented here precisely to demonstrate the core thesis of this paper: that meaningful progress must be local to deal with the interconnected challenges described earlier. This case study illustrates how the principles of ethical and sustainable practice can be enacted even within the constraints of a Western, Global North institution, proving that there is no universal solution, only context-specific ones. Indeed, even the EiMP finding its grounding in CUEiMS proved to be useful to deal with local institutional constraints: at the University of Cambridge, where the project was initially hosted[4], student societies are allowed to host non-examinable and optional seminars, and to invite guest speakers, even when other institutional constraints restrict the proper integration of topics into the examinable curriculum. Thus, students did not just learn from the project, but actively protected it.

We focus on the EiMP specifically because transforming the dominant Western mathematics — the paradigm often exported globally — requires localised interventions within the institutions that hold significant power. If ethical and sustainable consciousness cannot be fostered within these centers of influence, the broader project of mathematical practice focused on ethical and sustainable futures is likely to remain marginalised.

## The project's intellectual foundations

In the following, we begin to outline how the project developed a multi-layered perspective suitable for their local context, before we discuss what happened in an advanced topics seminar teaching the material.

At the macro level, addressing the challenge posed by entrenched philosophical commitments to mathematical neutrality, the project created the *Manifesto for the Responsible Development of Mathematical Works* (Chiodo & Müller, 2025). This document

---

[4] It was, from 2018 to 2024, known as the "Cambridge University Ethics in Mathematics Project": https://web.archive.org/web/20241001041556/http://www.ethics.maths.cam.ac.uk/

— "draw[ing] [among others] inspiration from more technology-specific ethics frameworks" (ibid., p. 13) — serves as a guideline for ethical mathematical practice. However, rather than proposing an objective-driven taxonomy, it delineates an integrated process featuring "10 pillars for responsible development". These 10 pillars lead mathematicians throughout the complete trajectory of mathematical work, starting with "Deciding whether to begin: Why are you providing this mathematical product or service, and should you even do so?" (ibid., p. 15), and ending with emergency response strategies for when things have gone wrong. The manifesto directly contests conceptions of value-neutrality by asserting that mathematical artefacts possess inherent political dimensions (pillar 9) and that mathematical practitioners bear responsibility for the ramifications of their work. The list of 10 pillars, around which the manifesto is built, are as follows:

1. Deciding whether to begin
2. Diversity and perspectives
3. Handling data and information
4. Data manipulation and inference
5. The mathematisation of the problem
6. Communicating and documenting your work
7. Falsifiability and feedback loops
8. Explainable and safe mathematics
9. Mathematical artefacts have politics
10. Emergency response strategies

Within the project and its teaching, the manifesto is used to demonstrate that justice-oriented and sustainable mathematical practices are not merely aspirational, but also systematically operationalisable within mainstream forms of mathematics, thereby advocating to reconceptualise the subject's philosophical foundations by taking the first step. As an iterative process-driven framework, it tries to challenge mathematicians to reflexively interrogate their problem selection criteria, the assumptions embedded within their mathematical solutions, and the normative commitments their work advances; all aimed at pushing the community further towards professional accountability. Importantly, the insights leading to the manifesto are not just practical experiences, but also draw on the argument that calls for a mere Hippocratic Oath for mathematicians and educators may be insufficient as the international mathematical communities lack the institutional infrastructure and ethical awareness which medicine had built over time (Müller et al., 2022; Rittberg, 2023b). Thus, behind the manifesto also lies the idea that internal critiques are necessary when external critiques are missing, but that they cannot fully replace them, and, indeed, should ideally go hand in hand.

At the meso level, the project created the *Four Levels of Ethical Engagement* framework (Chiodo & Bursill-Hall, 2018). This framework provides a conceptual architecture enabling students and professionals to comprehend and categorise different states of social and ethical responsibility. Transcending a binary classification of mathematicians' engagement with ethics as "ethical/unethical", it proposes a range from "Level 0" (ethical denialism in mathematics) to "Level 4" (proactive intervention against mathematical misappropriation). The levels are as follows:

- Level 0: Believing there is no ethics in mathematics
- Level 1: Realising there are ethical issues inherent in mathematics
- Level 2: Doing something: speaking out to other mathematicians
- Level 3: Taking a seat at the tables of power
- Level 4: Calling out the bad mathematics of others

This establishes a shared set of discursive conventions when other oaths or rules may be lacking in such guidance (e.g., Müller et al., 2022; Rycroft-Smith et al., 2024). By providing mathematicians and students with an explicit description of one developmental pathway for forming ethical consciousness, it is intended to enable them to progress from passive awareness to active engagement. The goal is to provide people, who will have only ever experienced mathematics presented as neutral, rational and objective at school and university, with a *positive* image of mathematics as a value-laden practice. While these levels describe where and how ethical (activism) can occur, they are not prescriptive, nor does the framework argue completeness. As Rycroft-Smith et al. (2024) show when bringing these levels into school settings, they have to be adjusted (and expanded) to fit other local circumstances.

At the micro level, fostering the ethical transformation of students while navigating institutional curricular constraints, the project created the pedagogical strategy of embedding ethical considerations within conventional mathematical exercises (Chiodo & Bursill-Hall, 2019; Chiodo et al., 2025). Recognising that separate ethics courses risk marginalising ethics and sustainability, this approach argues to integrate ethical questions, sustainability concerns, and sociocultural contexts directly into the standard mathematics courses. This interweaving of concerns aims to acclimate students in recognising and analysing the ethical and sustainable issues as constitutive elements of mathematical practice rather than as post-hoc considerations. In other words, the goal is to habitualise ethics and sustainability within existing institutional constraints (Chiodo & Bursill-Hall, 2019).

At the centre of this is Winner's (1980) insight that technical artefacts have politics. It is not just the mathematical practice which requires ethical and political decision-making, but that such decisions will lead to embedded politics in the mathematical artefact itself (Müller & Chiodo, 2023). The aim of the exercises is not just to teach students to move away from universal judgements, but also that good mathematical solutions are inherently dependent on local contexts. This includes designing exercises which teach students about different forms of aesthetic judgements within mathematics, environmental circumstances, different professional contexts, etc. These exercises can lead to what Shah (2025) calls "Teaching sustainability in mathematics problems? You must be joking!", whereby the exercises are misunderstood by mathematicians until they see how the mathematical component is strengthened by the ethical and sustainable considerations embedded within the questions and solutions.

Finally, to support individual students and mathematicians in navigating this complex terrain, the project produced a *Field Guide to Ethics in Mathematics* (Chiodo & Müller, 2024a). This resource aims to mediate between the diverse stakeholder interests and epistemological perspectives within different mathematical communities. It functions as a disciplinary primer, demonstrating that ethics in mathematics constitutes a "tractable" and "well-defined" domain of inquiry. By presenting different viewpoints and domains of concern, the field guide highlights the existence and legitimacy of diverse perspectives. The guide explores how questions surrounding ethics in mathematics typically fall into three concerns, or a combination thereof: concerns about mathematical practice and knowledge, community concerns, and concerns about the wider society and our planet. Thus, it presents ethics in mathematics as something fully localisable: it is not just that everyone can find something that is of concern to them, but that their local concerns are valid, and ought to be respected within different international mathematical communities, and that mathematical practice can be adjusted to fit local circumstances.

## The first time teaching the manifesto

Over time, and with EiMP's first release of the Manifesto in 2023, it was realised that the knowledge base of ethics in mathematics had developed so far, and become so operationalised, since the initial seminar series hosted by CUEiMS that an additional "advanced" series was set up. The EiMP found that it could not only teach students about ethical conduct as a mathematician, but also *how* to enact it. And so a non-examinable and optional advanced seminar series was formed, taking a physical-world case study of mathematical work and leading the students through the 10 pillars of the manifesto based on this. In its first year, the advanced series took the example of developing an AI-powered bus route and timetabling system (BRTS), and looked at how one would carry out such a project responsibly. Two of the students wrote a publicly available blog post about their experience (Yasmine & Siewert, 2024), from which we now analyse selected sections. The students were initially surprised that designing a BRTS can lead to complex ethical concerns:

> "At first glance, we did not appreciate how even something as simple as a BRTS could have various ramifications for citizens, democracy, and the environment. Yet, these rich ethical dilemmas occupied over 15 hours of discussions between students of maths, policy, law, and philosophy."

Their surprise illustrates again how traditional mathematics education can render even justice-centred topics such as public transport as supposedly neutral, and therefore mathematics *aimed* at the public as supposedly free of risk or potential harm. The students were accustomed to mathematical tasks which strip away all context, rendering ethical, social, and ecological consequences invisible. The 10 pillars successfully dismantled this by revealing the "wicked" nature of a BRTS. The fact that some students dedicated 15 hours to an optional ethics seminar, for which they did not receive any credits within their normal curriculum, demonstrates that the framework resonated strongly with their lived experiences and interests, prompting deep, interdisciplinary engagement. In effect, this displayed that the manifesto can act as a teaching framework focused on the existential sustainability concerns of a group of mathematics students. Writing about the first session, they note:

> "When presented with the task to develop a new BRTS, we were tempted to jump straight into finding solutions. However, Chiodo and Müller's pillars outline some practical issues that should be considered throughout the entire lifecycle of mathematical development. For instance, the first pillar, 'Deciding Whether to Begin', prompted us to do something which was very unfamiliar to the mathematicians in our group: scrutinising whether the question itself was a problem that needed to be solved."

The students quickly learned to accommodate fundamental shifts in their mathematical practice. The "temptation to jump straight into finding solutions" reflects the default mode of mathematics students trained to solve exercises and homework problems in an unquestioning way, never scrutinising the origins, framing, or requests made within the question; an implicit belief that the question cannot be wrong. The first pillar actively disrupted this impulse, and forced the students to perform an "unfamiliar" act by questioning the necessity of the mathematical task itself. This experience demonstrates the possibility to teach higher mathematics in ways that prioritise critical judgment and responsible problem-framing over mere abstraction and technical execution. This shift is more than just a change in procedure; it fosters a different professional identity, where the mathematician is a critical interrogator of problems rather than merely a solver of them. Operationalising ethical and sustainable thinking right at the inception of a given task is also necessary to appreciate the quest for diversity represented by pillar 2:

> "Next, for pillar two, 'Diversity and Perspectives', we looked inward at our own team's biases and recognised that, for instance, none of us had any experience in city planning, none of us knew the challenges of using the bus as a visually impaired person, nor had any of us ever driven a bus. Thus, we devised basic strategies on how we might reach out to others to ensure we do not overlook the needs of others or fundamental realities of designing public transport."

By applying the second pillar, the students quickly developed epistemic humility. They recognised the limitations of their own perspectives, and attempted to actively consider the needs of diverse stakeholders, including those groups who may be dependent on public transport, such as the visually impaired, allowing them to build a strong foundation of the principles of data ethics represented in pillars 3 and 4. What impacted the students quite remarkably was the idea that they might need to consult a bus driver in such work; someone who (likely) had far less mathematical knowledge than them, but could nonetheless impart great insight on how to view and approach the problem. They had not considered that a bus driver could help with their mathematics. Additionally, pillar 2 already connected the students with their out-of-university knowledge: as neither of them was visually impaired nor knew what it meant to use a bus as a visually impaired person, they had to fall back onto other, non-mathematical forms of knowledge and experiences. This perspective continued with the data-centric pillars 3 and for, about which they note:

> "[We] quickly found ourselves in conversation about whether we had legal access to use the data, whether our data distribution matched our application domain, and how our data could be used in harmful ways – for instance, to exclude certain members of the public from using the bus."

The epistemic humility and diversity-consciousness established in pillar 2 was maintained, and the "needs of others" explored previously were now used to help to include, rather than exclude, people from public transport. The focus was thus on the needs of the local community using the buses. Overall, through pillars 1-4 they started to see that shortcomings and oversights could creep into their work before they even carried out their first calculation or other mathematical step. Leading into the more mathematical pillar 5, they write:

> "So far [...] almost no mathematics has been done, which just shows one of the reasons ethical concerns are sometimes undervalued, namely that discussing them takes a lot of time which may be construed as "not work". The mathematics proper then begins in pillar five. 'The Mathematisation of the Problem' is about the ethics of the modelling and thus the mathematical methods used. For this point, it is most readily apparent that a huge amount of ethical challenges needs to tackled, for example, whether in our model there were any unjustifiable assumptions or reductions, how we approached nuances that could not be measured or quantified easily with a formula, and even details about the tools used such as computational cost, different degrees of confidence in models, and concerns regarding maintenance."

Privileging mathematical over ethical and sustainability considerations constitutes a substantial barrier within the subject's disciplinary culture, wherein temporal investments in such critical reflection may be marginalised as falling outside the bounds of legitimate academic labour. Nevertheless, students quickly recognised that the modelling process requires a deep commitment to its ethical and sustainable dimensions, seeing that one can carry out mathematical processes and get "nice" answers that are actually not at all fit for purpose. This was the point where the students saw that mathematics itself was the subservient actor to the actual problem at hand, rather than seeing the initial problem as a mere motivation to do some interesting mathematics. The students' language resonates with calls to revive the "cultural commons," which are the intergenerational, non-monetised skills and forms knowledge within a community focused on self-reliance. Bowers (2017 p. 54), notes that the "[r]evitalization of the cultural commons also requires recognizing the many ways capitalism, and its guiding ideology, attempts to integrate them into the market system." A crucial prerequisite for the long-term thinking required in the second half of the manifesto, focused on the the ethics of the mathematical post-production process:

> "Pillars seven through nine are about some key recurring issues which are not inherently mathematical, but inherent in all mathematical development and need to be kept in mind. By now we had quite some experience and were able to spend an hour listing issues around falsifiability [...], feedback loops [...], explainability [...], safety [...], and politics [...]"

Applying the second half of the framework gave the students lasting competence and confidence. They had internalised the process, and began to proactively consider complex issues, progressing from being guided by the framework to independently applying more substantive ethical reasoning. Once they had seen in the first half of the framework that there are some non- (or semi-) technical issues that need to be addressed *before* carrying out mathematical computations and modelling, it was almost trivial for them to see and accept that such issues also arose *after* the computations and modelling. Thus, these lessons underscored that ethical reasoning is a skill that can be learned and practiced. In other words, they levelled up in Chiodo & Bursill-Hall's (2018) level-based framework, and by the end were actively capable of calling out socially-unjust and unsustainable mathematics. About considering what to do when things have already gone wrong, the students hence wrote:

"In the final tenth pillar, 'Emergency Response Strategies', we again confronted an issue that is unfamiliar to mathematicians: "What could go wrong, and how can we respond?". A typical response to this is "I can prove that nothing will go wrong" or even the equivalent "we are clever, we will be careful to make sure nothing goes wrong". Therefore, we discussed some of the safeguards we might have in place to protect individuals if our BRTS fails and the appropriate response plans. This final pillar really focused our attention on finding non-mathematical interventions to engage when mathematics might fail."

This was probably one of the most counterintuitive steps for the students, requiring deep intellectual humility to accept that something could still go wrong, despite their best efforts. The mathematician's mantra of "we can do things perfectly" could no longer be retained at this point. In addition, the students saw the need to accept that additional mathematical work was no longer a useful action at this point, and that mathematics was not the solution to a mathematically-induced problem. With gentle reassurance, they came to terms with the fact that they could (and should) use mechanisms other than mathematics to address such issues; a broadening of perspective that they seldom receive in their standard mathematical training.

## Lessons learned

While, from the perspective of the authors, this multidimensional approach represents a necessary intervention within the Western mathematical research community, it is crucial to recognise it as a localised paradigm, tailored to succeed within a very specific set of circumstances. Born from practical concerns (Ernest, 2021b) and as a methodology "from mathematicians for mathematicians" (Müller, 2024, p. 83), its success hinges on a series of deliberate compromises. It offers critical-pragmatic tools that focus on individual responsibility and preventing harm from a technocratic standpoint; essentially presenting a version of sustainability and ethics digestible for mathematicians and students with little additional training outside of mathematics. These compromises were strategic, and seen as necessary to gain traction within a community historically resistant to philosophical and political critique.

This very pragmatism reveals its limitations when viewed from outside the context of a Western research university. The framework's philosophically light touch and its focus on the individual actions of powerful practitioners (Ernest, 2021b) may be insufficient for communities grappling with the deeper, systemic issues of coloniality, power, and epistemic injustice in mathematics (cf. Müller, 2024; Rittberg et al., 2020; Rittberg, 2023a, 2024). Thus, it presents *one* localised paradigm that may work for a research university in the Global North but it does not necessarily have to succeed in other places, where, as described earlier, powerful alternatives emerged due to their local sustainability and ethics concerns being different.

Even within universities of the Global North contexts can vary immensely. The local concerns that the project aims to address are seeing ethics, and bringing ethical awareness to a potentially unaware crowd in a practical fashion. The project attempts to equip high-performing mathematics students and professional mathematicians with a new professional identity that makes it easier for them to see other concerns. However, this still does not mean that these concerns match, or that the methodology created by the EiMP is necessarily applicable in other circumstances.

We also wish to recognise the significant burden that localised work places on educators. While localised, culturally-situated curricula are essential for ethical and sustainable mathematics education, they demand far more from its educators than traditional standardised approaches. As our experience within the EiMP highlights, mathematicians and educators must suddenly become curriculum designers, cultural researchers, and facilitators of complex ethical discussions. This presents individual mathematicians and educators with a substantial shift in focus, in new forms of knowledge, and in the demand to acquire new skills. This includes getting a deeper understanding of local contexts, the creation of custom teaching materials, and the ability to navigate sensitive topics within a mathematics classroom or lecture.

Finally, we acknowledge that a risk of this critical-pragmatic strategy is that it may reinforce the existing power structures by making the dominant mathematical paradigm more palatable. We see a tension between reform and revolution, and acknowledge that in the history of mathematics, revolutions can be difficult and rare (cf. Gillies, 1992). Our experience is that radical critique is always also relative to local contexts: what appears radical to a pure mathematician, might not appear radical to an educator focused on social justice. We operated on the premise that to enact change in students who dream about becoming mathematicians, and those who have become it, requires finding a tricky balance between activist ideals and pragmatic realities.

## Conclusion

We have shown that the urgent need to integrate ethics and sustainability into mathematics and its education is confronted by at least four deeply interconnected systemic challenges. These include the entrenched philosophical belief in mathematical neutrality which fosters resistance from practitioners, the difficulty of simultaneously reforming both the discipline of mathematics and its pedagogy, the persistent gap between academic theory and classroom practice that results in fragmented, patchwork implementations, and the profound imperative for epistemic decolonisation to avoid reproducing colonial epistemologies. We have argued that global frameworks like the United Nations' Sustainable Development Goals are ultimately insufficient for this, as they are rooted in a Western-centric development paradigm that can reproduce colonial hierarchies and fail to resolve inherent contradictions, such as the conflict between universal economic growth and ecological integrity; potentially only prolonging unsustainable and unethical modern conceptions of mathematics.

In contrast, localised interventions can offer a critical and pragmatic, multi-level strategy for cultivating ethical consciousness within specific communities. We discussed the Ethics in Mathematics Project as one such approach. In doing so, we noted that because its tools are designed to be operational and gain traction among mathematicians resistant to politicisation, this very pragmatism limits its capacity to address issues outside of specialised settings of higher mathematics. Ultimately, we argued that the path forward demands a departure from universalist, one-size-fits-all solutions, and to follow the call of ethnomathematics to embrace the diversity of culturally-situated, localised mathematical practices.

## Statements and Declarations

**AI Usage statement:** The authors used artificial intelligence to support the preparation of this manuscript. Specifically, Elicit and Gemini Deep Research were used to assist with the literature search. For improving grammar, clarity, and language, the authors further employed Grammarly, Gemini, and Claude. The authors reviewed and edited all AI-assisted content and take full responsibility for the final text.

**Acknowledgements:** The authors wish to thank their previous students and members of the Cambridge University Ethics in Mathematics Society for their continued support.

**Sources of funding:** The authors received no financial support to conduct this research.

**Financial or non-financial interests:** The authors have no financial or non-financial interests to declare.

**Conflicts of Interest/Competing interests:** The authors have no competing interests to declare.

**Ethical approval:** Not applicable.

**Availability of data and material:** Not applicable.

**Code availability:** Not applicable.

**Consent to participate:** Permission was obtained to study the students' blog post, and write about their experience, as part of this paper.

**Consent for publication:** Not applicable.